\documentclass[a4paper,11pt,english]{smfart}

\usepackage{amssymb}
\usepackage{url}
\usepackage[T1]{fontenc}
\usepackage[cal=boondoxo]{mathalfa}

\usepackage{fancybox}
\usepackage{slashbox,booktabs,amsmath}
\usepackage{makecell}

\usepackage[leftbars]{changebar}
\usepackage{todonotes}

\usepackage{palatino}
\usepackage{rotating}
\usepackage{graphicx}

\usepackage{enumerate}
\usepackage{enumitem}

\usepackage{color}
\definecolor{shadecolor}{gray}{0.90}

\def\bfit{\bfseries\itshape}

\def\proofname{D\'emonstration}

\input xypic
\xyoption{all}
\xyoption{arc}

\makeindex

\def\indexname{Index des notations}

\newtheorem{theo}{Theorem}[section]
\def\thetheo{\thesection.\arabic{theo}}
\newtheorem{prop}[theo]{Proposition}

\newtheorem{lem}[theo]{Lemma}

\newtheorem{coro}[theo]{Corollary}

\newtheorem{conj}[theo]{Conjecture}

\renewcommand\theequation{\thetheo}
\def\equat{\refstepcounter{theo}\begin{equation}}
\def\endequat{\end{equation}}

\renewcommand\thetable{\thetheo}

\renewcommand\thesection{\arabic{section}}
\renewcommand\thesubsection{\thesection.\Alph{subsection}}

\def\contentsname{Table des mati\`eres}
\def\listfigurename{Liste des figures}
\def\listtablename{Liste des tables}
\def\appendixname{Appendice}

\def\refname{R\'ef\'erences}

    \def\CM{{\mathbb{C}}}

\def\FG{{\mathfrak F}}

  \def\pG{{\mathfrak p}}

\def\SG{{\mathfrak S}}

    \def\ZM{{\mathbb{Z}}}

    \def\AC{{\mathcal{A}}}
    
    \def\CC{{\mathcal{C}}}
    \def\DC{{\mathcal{D}}}
    \def\EC{{\mathcal{E}}}
    
\def\Gb{{\mathbf G}}    
\def\Hb{{\mathbf H}}    
    \def\IC{{\mathcal{I}}}
    \def\JC{{\mathcal{J}}}
    
\def\Lb{{\mathbf L}}

\def\Pb{{\mathbf P}}    \def\PC{{\mathcal{P}}}

    \def\UC{{\mathcal{U}}}

    \def\XC{{\mathcal{X}}}
    
\def\Zb{{\mathbf Z}}    \def\ZC{{\mathcal{Z}}}

\def\Nrm{{\mathrm{N}}}

\def\Srm{{\mathrm{S}}}    
\def\Trm{{\mathrm{T}}}

\def\Zrm{{\mathrm{Z}}}

\def\a{\alpha}

\def\G{\Gamma}

\def\D{\Delta}
\def\e{\varepsilon}
\def\ph{\varphi}
\def\l{\lambda}
\def\L{\Lambda}

\def\o{\omega}
\def\O{\Omega}

\def\t{\tau}

\def\z{\zeta}

          \def\chit{{\tilde{\chi}}}

\def\mub{{\boldsymbol{\mu}}}

\def\Omeb{{\boldsymbol{\Omega}}}

\DeclareMathOperator{\Id}{{\mathrm{Id}}}

\DeclareMathOperator{\Irr}{{\mathrm{Irr}}}
\DeclareMathOperator{\Ker}{{\mathrm{Ker}}}

\DeclareMathOperator{\Part}{{\mathrm{Part}}}

\DeclareMathOperator{\Res}{{\mathrm{Res}}}

\DeclareMathOperator{\Tr}{{\mathrm{Tr}}}

\def\to{\rightarrow}
\def\longto{\longrightarrow}
\def\injto{\hookrightarrow}

\def\vide{\varnothing}
\def\DS{\displaystyle}

\def\SSS{\scriptscriptstyle}

\def\finl{~$\blacksquare$}

\def\lexp#1#2{\kern\scriptspace\vphantom{#2}^{#1}\kern-\scriptspace#2}
\def\le{\hspace{0.1em}\mathop{\leqslant}\nolimits\hspace{0.1em}}
\def\ge{\hspace{0.1em}\mathop{\geqslant}\nolimits\hspace{0.1em}}

\mathchardef\inferieur="321E
\mathchardef\superieur="321F

\def\eqna{\begin{eqnarray*}}
\def\endeqna{\end{eqnarray*}}

\def\maxi{{\mathrm{max}}}

\def\itemth#1{\item[${\mathrm{(#1)}}$]}

\catcode`\@=11
\long\def\@car#1#2\@nil{#1}
\long\def\@first#1#2{#1}
\long\def\@second#1#2{#2}
\long\def\ifempty#1{\expandafter\ifx\@car#1@\@nil @\@empty
  \expandafter\@first\else\expandafter\@second\fi}
\catcode`\@=12

\renewcommand{\Ref}{{\mathrm{Ref}}}

\def\boitegrise#1#2{\begin{centerline}{\fcolorbox{black}{shadecolor}{~
    \begin{minipage}[t]{#2}{\vphantom{~}#1\vphantom{$A_{\DS{A_A}}$}}
            \end{minipage}~}}\end{centerline}\medskip}

\def\ve{{\SSS{\vee}}}

\theoremstyle{remark}
\newtheorem{rema}[theo]{Remark}

\newtheorem{exemple}[theo]{Example}

\theoremstyle{plain}

\def\BIL{LR}
\def\GAUCHE{L}
\def\CAR{CAR}
\def\FAM{FAM}

\def\reg{{\mathrm{reg}}}

\def\xyinj{\ar@{^{(}->}}
\def\xysur{\ar@{->>}}

\bigskip

\makeatletter
\def\hlinewd#1{%
\noalign{\ifnum0=`}\fi\hrule \@height #1 %
\futurelet\reserved@a\@xhline}
\makeatother

\newlength\epaisLigne

\usepackage{multirow}

\newcommand{\longiso}{\stackrel{\sim}{\longrightarrow}}

\def\attractif{{\mathrm{att}}}

\makeatletter
\def\hlinewd#1{%
\noalign{\ifnum0=`}\fi\hrule \@height #1 %
\futurelet\reserved@a\@xhline}
\makeatother

\usepackage{multirow}

\usepackage{lscape}

\def\la{\langle}
\def\ra{\rangle}

\def\codim{\operatorname{codim}\nolimits}

\def\maxi{{\mathrm{max}}}

\begin{document}


\title{Regular automorphisms and Calogero-Moser families}

\author{{\sc C\'edric Bonnaf\'e}}
\address{IMAG, Universit\'e de Montpellier, CNRS, Montpellier, France} 

\makeatletter
\email{cedric.bonnafe@umontpellier.fr}
\makeatother

\date{\today}

\thanks{The author is partly supported by the ANR: 
Projects No ANR-16-CE40-0010-01 (GeRepMod) and ANR-18-CE40-0024-02 (CATORE).}

\pagestyle{myheadings}

\markboth{\sc C. Bonnaf\'e}{\sc Regular automorphisms and Calogero-Moser families}


%

\begin{abstract}
We study the subvariety of fixed points of an automorphism of a 
Calogero-Moser space induced by a regular element of finite order of 
the normalizer of the associated complex reflection group $W$. 
We determine some of (and conjecturally all) 
the ${\mathbb{C}}^\times$-fixed points of its unique irreducible component 
of maximal dimension in terms of the character table of $W$. 
This is inspired by the mysterious relations between 
the geometry of Calogero-Moser spaces and unipotent representations 
of finite reductive groups, which is the theme of another paper~\cite{cm-unip}. 
\end{abstract}

\maketitle

\bigskip

Let $V$ be a finite dimensional vector space and let $W$ be a finite 
subgroup of $\Gb\Lb_\CM(V)$ generated by reflections. To some parameter $k$, 
Etingof and Ginzburg~\cite{EG} have associated a normal 
irreducible affine complex variety $\ZC_k=\ZC_k(V,W)$ called 
a (generalized) {\it Calogero-Moser space}. If $\t$ is an 
element of finite order of the normalizer of $W$ in $\Gb\Lb_\CM(V)$ 
stabilizing the parameter $k$, it induces an automorphism 
of $\ZC_k$. 

We denote by $V_\reg$ 
the open subset of $V$ on which $W$ acts freely, and we assume that 
$V_\reg^\t \neq \vide$ (then $\t$ is called {\it regular}). In this case, there exists a unique 
irreducible component $(\ZC_k^\t)_\maxi$ of $\ZC_k^\t$ of maximal 
dimension (as it will be 
explained in Section~\ref{sec:regular}). Recall that $\ZC_k$ 
is endowed with a $\CM^\times$-action and that we have 
a surjective map $\Irr(W) \to \ZC_k^{\CM^\times}$ defined by Gordon~\cite{gordon} 
(induced by the action of the center of a rational Cherednik algebra on {\it baby Verma modules}) 
whose fibers are called the {\it Calogero-Moser $k$-families} of $W$. 
If $p \in \ZC_k^{\CM^\times}$, we denote by $\FG_p$ its associated 
Calogero-Moser $k$-family. 
It is a natural question to wonder which $\CM^\times$-fixed points 
of $\ZC_k^\t$ belong to $(\ZC_k^\t)_\maxi$. The aim of this note 
is to provide a partial answer in terms of the character table of $W$:

\bigskip

\noindent{\bfit Theorem A.~---} {\it Assume that $V_\reg^\t \neq \vide$. 
Let $p \in \ZC_k^{\CM^\times}$ be such that $\t(p)=p$. 
If $\sum_{\chi \in \FG_p^\t} |\chit(\t)|^2 \neq 0$, then 
$p \in (\ZC_k^\t)_\maxi$.}

\bigskip

In this statement, if $\chi$ is a $\t$-stable irreducible character 
of $W$, we denote by $\chit$ an extension of $\chi$ to the finite 
group $W\langle \t \rangle$ (note that $|\chit(\t)|^2$ does 
not depend on the choice of $\chit$). Our proof of Theorem~A makes an extensive use 
of the Gaudin operators introduced in~\cite[\S{8.3.B}]{calogero}. 
This result is also inspired by the theory of unipotent 
representations of finite reductive groups and some 
conjectures of Brou\'e-Michel~\cite{broue-michel} on the cohomology 
of Deligne-Lusztig varieties associated with regular elements 
in the sense of Springer~\cite{springer} and by~\cite[Rem.~4.21]{spetses} 
(this will also be discussed in~\cite{cm-unip}). If we believe 
in this analogy, we can conjecture that the converse of Theorem~A holds:

\bigskip

\centerline{\begin{minipage}{0.7\textwidth}
\noindent{\bfit Conjecture B.~---} {\it Assume that $V_\reg^\t \neq \vide$. Let $p \in \ZC_k^{\CM^\times}$ 
be such that $\t(p)=p$. Then $p \in (\ZC_k^\t)_\maxi$ 
if and only if $\sum_{\chi \in \FG_p^\t} |\chit(\t)|^2 \neq 0$.}
\end{minipage}}

\bigskip

\noindent{\bf General notation.} 
Throughout this paper, we will abbreviate $\otimes_\CM$ as 
$\otimes$ and all varieties will be algebraic, complex, quasi-projective 
and reduced. If $\XC$ is an affine variety, we denote by $\CM[\XC]$ 
its coordinate ring.

If $X$ is a subset of a vector space $V$ (or of its dual $V^*$), and if $\G$ is a subgroup 
of $\Gb\Lb_\CM(V)$, we denote by $\G_X$ the pointwise stabilizer of $X$. 
If moreover $\G$ is finite, we will identify $(V^\G)^*$ and $(V^*)^\G$. 

\bigskip

%

\section{Set-up}\label{sec:notation}

\medskip

\boitegrise{{\bf Hypothesis and notation.} {\it We fix in this paper a finite dimensional 
complex vector space $V$ and a finite subgroup $W$ of $\Gb\Lb_\CM(V)$. We set
$$\Ref(W)=\{s \in W~|~\codim_\CM V^s=1\}$$
and we assume throughout this paper that
$$W=\langle \Ref(W) \rangle,$$
i.e. that $W$ is a complex reflection group.}}{0.75\textwidth}

\bigskip

\subsection{About ${\boldsymbol{W}}$} 
We set $\e : W \to \CM^\times$, $w \mapsto \det(w)$. 
We identify $\CM[V]$ (resp. $\CM[V^*]$) with the symmetric 
algebra $\Srm(V^*)$ (resp. $\Srm(V)$).

We denote by $\AC$ the set of {\it reflecting hyperplanes} of $W$, namely
$$\AC=\{V^s~|~s \in \Ref(W)\}.$$
If $H \in \AC$, we denote by $\a_H$ an element of $V^*$ such that 
$H=\Ker(\a_H)$ and by $\a_H^\vee$ an element of $V$ such that 
$V=H \oplus \CM \a_H^\vee$ and the line $\CM\a_H^\vee$ is $W_H$-stable.
We set $e_H=|W_H|$. Note that $W_H$ is cyclic of order $e_H$ and that 
$\Irr(W_H)=\{\Res_{W_H}^W \e^j~|~0 \le j \le e-1\}$. We denote by $\e_{H,j}$ 
the (central) primitive idempotent of $\CM W_H$ associated with the character 
$\Res_{W_H}^W \e^{-j}$, namely
$$\e_{H,j}=\frac{1}{e_H}\sum_{w \in W_H} \e(w)^j w \in \CM W_H.$$
If $\O$ is a $W$-orbit of reflecting hyperplanes, we write $e_\O$ for the 
common value of all the $e_H$, where $H \in \O$. 
We denote by $\aleph$ the set of pairs $(\O,j)$ where $\O \in \AC$ and 
$0 \le j \le e_\O-1$. 
The vector space of families of complex numbers 
indexed by $\aleph$ will be denoted by $\CM^\aleph$: elements 
of $\CM^\aleph$ will be called {\it parameters}. 
If $k=(k_{\O,j})_{(\O,j) \in \aleph} \in \CM^\aleph$, we 
define $k_{H,j}$ for all $H \in \O$ and $j \in \ZM$ by 
$k_{H,j}=k_{\O,j_0}$ where $\O$ is the $W$-orbit of $H$ 
and $j_0$ is the unique element of $\{0,1,\dots,e_H-1\}$ such that 
$j \equiv j_0 \mod e_H$. 

\medskip

We denote by $V_\reg$ the set of elements $v$ of $V$ such that $W_v=1$. 
It is an open subset of $V$ and recall from Steinberg-Serre Theorem~\cite[Theo.~4.7]{broue} 
that 
\equat\label{eq:steinberg}
V_\reg=V \setminus \bigcup_{H \in \AC} H.
\endequat

\subsection{Rational Cherednik algebra at ${\boldsymbol{t=0}}$}
Let $k \in \CM^\aleph$. 
We define the {\it rational Cherednik algebra $\Hb_k$} ({\it at $t=0$}) to be the quotient 
of the algebra $\Trm(V\oplus V^*)\rtimes W$ (the semi-direct product of the tensor algebra 
$\Trm(V \oplus V^*)$ with the group $W$) 
by the relations 
\equat\label{eq:rels}
\begin{cases}
[x,x']=[y,y']=0,\\
[y,x]=\DS{\sum_{H\in\mathcal{A}} \sum_{j=0}^{e_H-1}
e_H(k_{H,j}-k_{H,j+1}) 
\frac{\langle y,\a_H \rangle \cdot \langle \a_H^\ve,x\rangle}{\langle \a_H^\ve,\a_H\rangle} \e_{H,j}},
\end{cases}
\endequat
for all $x$, $x'\in V^*$, $y$, $y'\in V$. 
Here $\la\ ,\ \ra: V\times V^*\to\CM$ is the standard pairing. 
The first commutation relations imply that 
we have morphisms of algebras $\CM[V] \to \Hb_k$ and $\CM[V^*] \to \Hb_k$. 
Recall~\cite[Theo.~1.3]{EG} 
that we have an isomorphism of $\CM$-vector spaces 
\equat\label{eq:pbw}
\CM[V] \otimes \CM W \otimes \CM[V^*] \longiso \Hb_k
\endequat
induced by multiplication (this is the so-called {\it PBW-decomposition}). 

\medskip

\begin{rema}\label{rem:parametres particuliers}
Let $(l_\O)_{\O \in \AC/W}$ be a family of complex numbers and let 
$k' \in \CM^\aleph$ be defined by $k_{\O,j}'=k_{\O,j} + l_\O$. Then 
$\Hb_k=\Hb_{k'}$. This means that there is no restriction to generality 
if we consider for instance only 
parameters $k$ such that $k_{\O,0}=0$ for all $\O$, 
or only parameters $k$ such that $k_{\O,0}+k_{\O,1}+\cdots+k_{\O,e_\O-1}=0$ 
for all $\O$ (as in~\cite{calogero}).\finl
\end{rema}

\medskip

\subsection{Calogero-Moser space} 
We denote by $\Zb_k$ the center of the algebra $\Hb_k$: it is well-known~\cite[Theo~3.3~and~Lem.~3.5]{EG} that 
$\Zb_k$ is an integral domain, which is integrally closed. Moreover, it contains 
$\CM[V]^W$ and $\CM[V^*]^W$ as subalgebras~\cite[Prop.~3.6]{gordon} 
(so it contains $\Pb=\CM[V]^W \otimes \CM[V^*]^W$), 
and it is a free $\Pb$-module of rank $|W|$. We denote by $\ZC_k$ the 
affine algebraic variety whose ring of regular functions $\CM[\ZC_k]$ is $\Zb_k$: 
this is the {\it Calogero-Moser space} associated with the datum $(V,W,k)$. 
It is irreducible and normal. 

We set $\PC=V/W \times V^*/W$, so that $\CM[\PC]=\Pb$ and the inclusion 
$\Pb \injto \Zb_k$ induces a morphism of varieties 
$$\Upsilon_k : \ZC_k \longto \PC$$
which is finite and flat.

\bigskip

\subsection{Calogero-Moser families}\label{sub:cm-families}
Using the PBW-decomposition, we define a $\CM$-linear map 
$\Omeb^{\Hb_k} : \Hb_k \longto \CM W$
by 
$$\Omeb^{\Hb_k}(f w g)=f(0)g(0)w$$
for all $f \in \CM[V]$, $g \in \CM[V^*]$ and $w \in \CM W$. This map is $W$-equivariant 
for the action on both sides by conjugation, so it induces a well-defined $\CM$-linear map 
$$\Omeb^k : \Zb_k \longto \Zrm(\CM W).$$
Recall from~\cite[Cor.~4.2.11]{calogero} that $\Omeb^k$ is a morphism of algebras. 

Calogero-Moser families were defined by Gordon 
using his theory of {\it baby Verma modules}~\cite[\S{4.2}~and~\S{5.4}]{gordon}. We explain  
here an equivalent definition given in~\cite[\S{7.2}]{calogero}. 
If $\chi \in \Irr(W)$, we denote by $\o_\chi : \Zrm(\CM W) \to \CM$ 
its central character (i.e., $\o_\chi(z)=\chi(z)/\chi(1)$ is the scalar by which 
$z$ acts on an irreducible representation affording the character $\chi$). 
We say that two characters $\chi$ and $\chi'$ belong to the same 
{\it Calogero-Moser $k$-family} 
if $\o_\chi \circ \Omeb^k = \o_{\chi'} \circ \Omeb^k$. 

In other words, the map $\o_\chi \circ \Omeb^k : \Zb_k \to \CM$ is a 
morphism of algebras, so it might be viewed as a point $\ph_k(\chi)$ of $\ZC_k$, 
which is easily checked to be $\CM^\times$-fixed. 
This defines a surjective map 
$$\ph_k : \Irr(W) \longto \ZC_k^{\CM^\times}$$
whose fibers are the Calogero-Moser $k$-families. 
If $p \in \ZC_k^{\CM^\times}$, we denote by $\FG_p$ 
the corresponding Calogero-Moser $k$-family.

\bigskip

\subsection{Other parameters} 
Let $\CC$ denote the space of maps $\Ref(W) \to \CM$ 
which are constant on conjugacy classes of reflections. 
The element
$$\sum_{(\O,j) \in \aleph} \sum_{H \in \O} (k_{H,j}-k_{H,j+1}) e_H \e_{H,j}$$
of $\Zrm(\CM W)$ is supported only by reflections, so there exists 
a unique map $c_k \in \CC$ such that 
$$\sum_{(\O,j) \in \aleph} \sum_{H \in \O} (k_{H,j}-k_{H,j+1}) e_H \e_{H,j}
= \sum_{s \in \Ref(W)} (\e(s)-1) c_k(s) s.$$
Then the map $\CM^\aleph \to \CC$, $k \mapsto c_k$ is linear and surjective. 
With this notation, we have 
\equat\label{eq:cs}
[y,x] = \sum_{s \in \Ref(W)} (\e(s)-1)\hskip1mm c_k(s) 
\hskip1mm
\frac{\langle y,\a_s \rangle \cdot \langle \a_s^\ve,x\rangle}{\langle 
\a_s^\ve,\a_s\rangle}
\hskip1mm s,
\endequat
for all $y \in V$ and $x \in V^*$. Here, $\a_s=\a_{V^s}$ and $\a_s^\ve=\a_{V^s}^\ve$.

\bigskip

\subsection{Actions on the Calogero-Moser space}
The Calogero-Moser space $\ZC_k$ is endowed with a $\CM^\times$-action and 
an action of the stabilizer of $k$ in $\Nrm_{\Gb\Lb_\CM(V)}(W)$, which 
are described below.

\bigskip

\subsubsection{Grading, $\CM^\times$-action}
The algebra $\Trm(V\oplus V^*)\rtimes W$ can be $\ZM$-graded in such a way that the 
generators have the following degrees
$$
\begin{cases}
\deg(y)=-1 & \text{if $y \in V$,}\\
\deg(x)=1 & \text{if $x \in V^*$,}\\
\deg(w)=0 & \text{if $w \in W$.}
\end{cases}
$$
This descends to a $\ZM$-grading on $\Hb_k$, because the defining relations~(\ref{eq:rels}) 
are homogeneous. Since the center of a graded algebra is always graded, the subalgebra $\Zb_k$ 
is also $\ZM$-graded.  So the Calogero-Moser space $\ZC_k$ 
inherits a regular $\CM^\times$-action. Note also that
by definition $\Pb=\CM[V]^W \otimes \CM[V^*]^W$ is clearly a graded 
subalgebra of $\Zb_k$. 

\bigskip

\subsubsection{Action of the normalizer}\label{sub:normalisateur}
The group $\Nrm_{\Gb\Lb_\CM(V)}(W)$ acts on the set $\aleph$ and so on 
the space of parameters $\CM^\aleph$. If $\t \in \Nrm_{\Gb\Lb_\CM(V)}(W)$, 
then $\t$ induces an isomorphism of algebras $\Hb_k \longto \Hb_{\t(k)}$. 
So, if $\t(k)=k$, then it induces an action on the algebra $\Hb_k$ 
(and so on its center $\Zb_k$ and on the Calogero-Moser space $\ZC_k$). 

We say that $\t$ is a {\it regular} element of $\Nrm_{\Gb\Lb_\CM(V)}(W)$ 
if $V_\reg^\t \neq \vide$. 

\bigskip

\boitegrise{{\bf Notation.} {\it From now on, and until the end of this paper, we 
fix a parameter $k \in \CM^\aleph$ and a {\bfit regular} element $\t$ of {\bfit finite} 
order of $\Nrm_{\Gb\Lb_\CM(V)}(W)$ such that $\t(k)=k$.}}{0.82\textwidth}

\bigskip

We denote by $\ZC_k^\t$ the variety of fixed points of $\t$ in $\ZC_k$, 
endowed with its reduced structure. All the above constructions are $\t$-equivariant: 
for instance, the map $\ph_k : \Irr(W) \longto \ZC_k^{\CM^\times}$ 
is $\t$-equivariant. 

Let us recall the following consequence~\cite[Prop.~3.5~and~Theo.~4.2]{springer} 
of the above hypothesis:

\bigskip

\begin{theo}[Springer]\label{theo:springer}
The group $W^\t$ acts as a reflection 
group on $V^\t$ and the natural map $V^\t/W^\t \to (V/W)^\t$ 
is an isomorphism of varieties.
\end{theo}

\bigskip

\begin{coro}\label{coro:springer}
The natural map $(V_\reg^\t \times V^{*\t})/W^\t \to ((V_\reg \times V^*)/W)^\t$ 
is an isomorphism of varieties.
\end{coro}

\medskip

\begin{proof}
Since $W$ acts freely on $V_\reg \times V^*$, the quotient $(V_\reg \times V^*)/W$ is smooth. 
Consequently, the variety of fixed points $((V_\reg \times V^*)/W)^\t$ is also smooth. 
Similarly, $(V_\reg^\t \times V^{*\t})/W^\t$ is smooth. Since a bijective 
morphism between smooth varieties is an isomorphism, we only need to show that 
the above natural map is bijective.

First, if $(v_1,v_1^*)$ and $(v_2,v_2^*)$ are two elements of $V_\reg^\t \times V^{*\t}$ 
belonging to the same $W$-orbit, there exists $w \in W$ such that 
$v_2=w(v_1)$. Since $v_1$ and $v_2$ are $\t$-stable, we also have 
$\t(w)(v_1)=v_2$, and so $v_1=w^{-1}\t(w)(v_1)$. Since $v_1 \in V_\reg$, this forces 
$\t(w)=w$ and the injectivity follows.

Now, if $(v,v^*) \in V_\reg \times V^*$ is such that its $W$-orbit is $\t$-stable, 
then the $W$-orbit of $v$ is $\t$-stable. So Theorem~\ref{theo:springer} shows 
that we may assume that $\t(v)=v$. The hypothesis implies that there exists 
$w \in W$ such that $\t(v)=w(v)$ and $\t(v^*)=w(v^*)$. But $\t(v)=v \in V_\reg$, 
so $w=1$. In particular, $\t(v^*)=v^*$, and the surjectivity follows.
\end{proof}

\section{Irreducible component of maximal dimension}\label{sec:regular}

\medskip

Let $(\ZC_k)_\reg$ denote the open subset $\Upsilon_k^{-1}(V_\reg/W \times V^*/W)$. 
By~\cite[Prop.~4.11]{EG}, 
we have a $\CM^\times$-equivariant and $\t$-equivariant isomorphism 
\equat\label{eq:ouvert}
(\ZC_k)_\reg \simeq (V_\reg \times V^*)/W.
\endequat
This shows that $(\ZC_k)_\reg$ is smooth and so $(\ZC_k)_\reg^\t$ is also smooth. 
By Corollary~\ref{coro:springer}, this implies that 
\equat\label{eq:zreg}
(\ZC_k)_\reg^\t \simeq (V_\reg^{\t} \times V^{*\t})/W^{\t}.
\endequat
In particular it is irreducible. We denote by $(\ZC_k^\t)_\maxi$ its closure: 
it is an irreducible closed subvariety of $\ZC_k^\t$. 

Moreover, $(\ZC_k)_\reg^\t$ has dimension $2\dim V^\t$ by Corollary~\ref{coro:springer}. 
So $\dim \ZC_k^\t \ge 2\dim V^\t = \dim (\ZC_k^\t)_\maxi$. But, on the other hand, 
$\Upsilon_k(\ZC_k^\t) \subset (V/W)^\t \times (V^*/W)^\t$. Since $\Upsilon_k$ 
is a finite morphism, we get from Theorem~\ref{theo:springer} that 
$\dim \ZC_k^\t \le 2 \dim V^\t$. Hence
\equat\label{eq:dim-zktau}
\dim \ZC_k^\t=\dim (\ZC_k^\t)_\maxi = 2 \dim V^\t.
\endequat
This shows that $(\ZC_k^\t)_\maxi$ is an irreducible component 
of maximal dimension of $\ZC_k^\t$. 

\bigskip

\begin{prop}\label{prop:zk-max}
The closed subvariety $(\ZC_k^\t)_\maxi$ of $\ZC_k^\t$ is the unique 
irreducible component of maximal dimension.
\end{prop}

\medskip

\begin{proof}
Let $\XC$ be an irreducible component of $\ZC_k^\t$ of dimension $2\dim V^\t$. 
Since $\Upsilon_k$ is finite, the image $\Upsilon_k(\XC)$ 
is closed in $V/W \times V^*/W$, irreducible of dimension 
$2 \dim(V^\t)$ and contained in $(V/W)^\t \times (V^*/W)^\t$. 
By Theorem~\ref{theo:springer}, we get that $\Upsilon_k(\XC)=(V/W)^\t \times (V^*/W)^\t$. 

Let $\UC=\Upsilon_k^{-1}(V_\reg/W \times V^*/W) \cap \XC$. Then $\UC$ is a non-empty open 
subset of $\XC$: since $\XC$ is irreducible, this forces $\UC$ to have 
dimension $2\dim(V^\t)$. But $\UC$ is contained in $(\ZC_k)_\reg^\t$ 
which is irreducible of the same dimension, so the closure of $\UC$ 
contains $(\ZC_k)_\reg^\t$. This proves that $\XC=(\ZC_k^\t)_\maxi$. 
\end{proof}

\bigskip

\begin{coro}\label{coro:zk-max}
$\Upsilon_k((\ZC_k^\t)_\maxi)=(V/W)^\t \times (V^*/W)^\t$.
\end{coro}

\bigskip

It is natural to ask which $\CM^\times$-fixed points of $\ZC_k$ belong 
to $(\ZC_k^\t)_\maxi$. Inspired by the representation theory of 
finite reductive groups (see~\cite{broue-michel} and~\cite[Rem.~4.21]{spetses}), 
we propose an answer to this question in terms of the character table 
of the finite group $W\langle \t\rangle$ (see~\cite[Ex.~12.9]{cm-unip} for 
some explanations). We first need some notation. 

If $\chi \in \Irr(W)$, we denote by $E_\chi$ a $\CM W$-module affording the 
character $\chi$. If moreover $\chi$ is $\t$-stable, we fix a structure 
of $\CM W\langle\t\rangle$-module on $E_\chi$ extending the structure of $\CM W$-module, 
and we denote by $\chit$ its associated irreducible character of $W\langle\t\rangle$. 
Note that the real number $|\chit(\t)|^2$ does not depend on the choice of $\chit$. 

\bigskip

\begin{quotation}
\begin{conj}\label{conj:chi-tau}
Recall that $\t$ is regular. 
Let $p \in \ZC_k^{\CM^\times}$ be such that $\t(p)=p$. Then $p$ belongs to $(\ZC_k^\t)_\maxi$ 
if and only if $\sum_{\chi \in \FG_p^\t} |\chit(\t)|^2  \neq 0$.
\end{conj}
\end{quotation}

\bigskip

\begin{rema}\label{rem:famille-tau}
Let $\FG$ be a $\t$-stable Calogero-Moser family. Then 
$\FG$ contains a unique irreducible character $\chi_\FG$ with minimal 
$b$-invariant~\cite[Theo.~7.4.1]{calogero}, where the $b$-invariant of an irreducible character 
$\chi$ is the minimal natural number $j$ such that $\chi$ occurs in the 
$j$-th symmetric power of the natural representation $V$ of $W$. 
From this characterization, we see that $\chi_\FG$ is $\t$-stable. 
In particular, any $\t$-stable Calogero-Moser family contains 
at least one $\t$-stable character.\finl
\end{rema}

\medskip

In general, we are only able to prove the ``if'' part of 
Conjecture~\ref{conj:chi-tau}.

\medskip

\begin{theo}\label{theo:chi-tau}
Recall that $\t$ is regular. 
Let $p \in \ZC_k^{\CM^\times}$ be such that $\t(p)=p$. If 
$\sum_{\chi \in \FG_p^\t} |\chit(\t)|^2  \neq 0$, then $p$ belongs to $(\ZC_k^\t)_\maxi$. 
\end{theo}

\medskip

The next two sections are devoted to the proof of Theorem~\ref{theo:chi-tau}. 

\bigskip

\section{Verma modules}

\medskip

\subsection{Definition}
Recall that $\CM[V] \rtimes W$ is a subalgebra of $\Hb_k$ 
(it is the image of $1 \otimes \CM W \otimes \CM[V]$ by the 
PBW-decomposition~\ref{eq:pbw}). If $E$ is a $\CM W$-module, 
we denote by $E^\#$ the $(\CM[V^*] \rtimes W)$-module 
extending $E$ by letting any element $f \in \CM[V^*]$ 
acting by multiplication by $f(0)$. 
If $\chi \in \Irr(W)$, we define an $\Hb_k$-module $\D(\chi)$ 
as follows:
$$\D(\chi)=\Hb_k \otimes_{\CM[V^*] \rtimes W} E_\chi^\#.$$
Then $\D(\chi)$ is called a {\it Verma module} of $\Hb_k$ 
(see~\cite[\S{5.4.A}]{calogero}: in this reference, $\D(\chi)$ 
is denoted by $\D(E_\chi^\#)$). Let $\Hb_k^\reg$ denote the localization 
of $\Hb_k$ at $\Pb_\reg=\CM[V_\reg/W] \otimes \CM[V^*]$. By~\cite[Prop.~4.11]{EG}, 
we have an isomorphism  
$\CM[V_\reg \times V^*]\rtimes W \simeq \Hb_k^\reg$. 
We denote by $\D^\reg(\chi)$ the localization of $\D(\chi)$ at 
$\Hb_k^\reg$. So, by restriction to 
$\CM[V_\reg \times V^*]$, the Verma module $\D(\chi)$ might be viewed as a 
$W$-equivariant coherent sheaf on $V_\reg \times V^*$. 
We also view $e\D(\chi)$ as a coherent sheaf on $\ZC_k$, so that 
$e\D^\reg(\chi)$ may be viewed as a coherent sheaf on $(V_\reg \times V^*)/W$. 
If $p \in \ZC_k$ (or if $(v,v^*) \in V_\reg \times V^*$), we denote by 
$e\D(\chi)_p$ (respectively 
$e\D(\chi)_{W\cdot(v,v^*)}=e\D^\reg(\chi)_{W \cdot (v,v^*)}$, 
respectively $\D^\reg(\chi)_{v,v^*}$) the restriction of $e\D(\chi)$ (respectively 
of $e\D(\chi)$ or $e\D^\reg(\chi)$, respectively $\D^\reg(\chi)$) 
at the point $p$ (respectively $W \cdot (v,v^*) \in (V_\reg \times V^*)/W \simeq (\ZC_k)_\reg$, 
respectively $(v,v^*)$). It follows from the definition 
that the support of $e\D(\chi)$ is contained in 
$\Upsilon_k^{-1}(V/W \times 0)$, and recall that, through the isomorphism 
$\ZC_k^\reg \simeq (V_\reg \times V^*)/W$, $\Upsilon_k^{-1}(V_\reg/W \times 0)$ 
is not necessarily contained in $(V_\reg \times \{0\})/W$. 

\bigskip

\begin{lem}\label{lem:fibre-fixe}
Let $\chi \in \Irr(W)$ and let $p \in \ZC_k^{\CM^\times}$. Then 
$e\D(\chi)_p \neq 0$ if and only if $\chi \in \FG_p$.
\end{lem}

\medskip

\begin{proof}
Let $\pG_0$ denote the maximal ideal of the algebra $\Pb=\CM[\PC]$ consisting 
of functions which vanish at $0$. Then $\D(\chi)/\pG_0 \D(\chi)$ 
is a representation of the {\it restricted rational Cherednik algebra} 
$\Hb_k/\pG_0\Hb_k$ which coincides with the {\it baby Verma module} 
defined by Gordon~\cite[\S{4.2}]{gordon}. As $\ZC_k^{\CM^\times}=\Upsilon_k^{-1}(0)$, 
the result follows from the very definition of Calogero-Moser families 
in terms of baby Verma modules and the 
fact that it is equivalent to the definition given in~\S\ref{sub:cm-families}.
\end{proof}

\subsection{Bialynicki-Birula decomposition}
We denote by $\ZC_k^\attractif$ the {\it attracting set} of $\ZC_k$ 
for the action of $\CM^\times$, namely
$$\ZC_k^\attractif=\{p \in \ZC_k~|~\lim_{\xi \to 0} \lexp{\xi}{p}~\text{exists}\}.$$
Recall from~\cite[Chap.~14]{calogero} the following facts:

\bigskip

\begin{prop}\label{prop:attractif}
With the above notation, we have:
\begin{itemize}
\itemth{a} The map $\lim : \ZC_k^\attractif \longto \ZC_k^{\CM^\times}$, 
$p \mapsto \lim_{\xi \to 0} \lexp{\xi}{p}$ is a morphism of varieties.

\itemth{b} $\ZC_k^\attractif=\Upsilon_k^{-1}(V/W \times \{0\})$.

\itemth{c} If $\IC$ is an irreducible component of $\ZC_k^\attractif$, 
then $\IC$ is $\CM^\times$-stable and $\Upsilon_k(\IC)=V/W \times \{0\}$ 
and $\lim(\IC)$ is a single point.

\itemth{d} If $\chi \in \Irr(W)$, then the support of $e\D(\chi)$ 
is a union of irreducible components of $\ZC_k^\attractif$.

\itemth{e} If $\IC$ is an irreducible component of $\ZC_k^\attractif$, then 
there exists $\chi \in \Irr(W)$ such that the support of $e\D(\chi)$ contains $\IC$. 
\end{itemize}
\end{prop}

\bigskip

We first propose a characterization of points $p \in \ZC_k^{\CM^\times}$ which belong 
to $(\ZC_k^\t)_\maxi$ in terms of Verma modules.

\bigskip

\begin{lem}\label{lem:zcmud}
Let $p \in \ZC_k^{\CM^\times}$ and assume that $\t(p)=p$. 
Then $p \in (\ZC_k^\t)_\maxi$ if and only if there exist 
$\chi \in \FG_p^k$ and $(v,v^*) \in V_\reg^\t \times V^{*\t}$ such that 
$e\D(\chi)_{W \cdot (v,v^*)} \neq 0$.
\end{lem}

\medskip

\begin{proof}
Let $(\ZC_k^\t)_\maxi^\attractif$ denote the attracting set of 
$(\ZC_k^\t)_\maxi$. Then Corollary~\ref{coro:zk-max} implies that 
$\Upsilon_k((\ZC_k^\t)_\maxi)=(V/W)^\t \times \{0\}$.  
Since $\Upsilon_k$ is a finite morphism, the same arguments used in~\cite[Chap.~14]{calogero} 
to prove the Proposition~\ref{prop:attractif} above yields the following 
statements:
\begin{itemize}
\itemth{a} The map $\lim : (\ZC_k^\t)_\maxi^\attractif \longto (\ZC_k^\t)_\maxi^{\CM^\times}$, 
$p \mapsto \lim_{\xi \to 0} \lexp{\xi}{p}$ is a morphism of varieties.

\itemth{b} $(\ZC_k^\t)_\maxi^\attractif=(\ZC_k^\t)_\maxi \cap 
\Upsilon_k^{-1}((V/W)^\t \times \{0\})$.

\itemth{c} If $\IC$ is an irreducible component of $(\ZC_k^\t)_\maxi^\attractif$, 
then $\IC$ is $\CM^\times$-stable and $\Upsilon_k(\IC)=(V/W)^\t \times \{0\}$ 
and $\lim(\IC)$ is a single point.
\end{itemize}

Assume that $p \in (\ZC_k^\t)_\maxi$. Let 
$\IC$ be an irreducible component of $(\ZC_k^\t)_\maxi^\attractif \cap \lim^{-1}(p)$. 
Then $\IC$ is contained in an irreducible component $\IC'$ of 
$(\ZC_k^\t)_\maxi^\attractif$. Since $\lim(\IC')$ is a single point by~(c), 
we have $\lim(\IC')=\{p\}$ and so $\IC=\IC'$. Still by~(c), this says that 
$\Upsilon_k(\IC)=(V/W)^\t \times \{0\}$. So let $q \in \IC$ be such that 
$\Upsilon_k(q) \in (V_\reg/W)^\t \times \{0\}$. 

Now, let $\JC$ be an irreducible component of $\ZC_k^\attractif$ containing 
$\IC$. By Proposition~\ref{prop:attractif}(e), there exists $\chi \in \Irr(W)$ 
such that the support of $e\D(\chi)$ contains $\JC$. In particular, 
$e\D(\chi)_p \neq 0$ and so $\chi \in \FG_p$ by Lemma~\ref{lem:fibre-fixe}. 
But also $e\D(\chi)_q \neq 0$. 
Since $q \in (\ZC_k^\t)_\maxi$ and $\Upsilon_k(q) \in V_\reg/W$, 
it follows 
that there exists $(v,v^*) \in V_\reg^\t \times V^{*\t}$ such that 
$e\D(\chi)_{W \cdot (v,v^*)} \neq 0$, as desired.

\medskip

Conversely, assume that there exist both $\chi \in \FG_p^k$ and 
$(v,v^*) \in V_\reg^\t \times V^{*\t}$ such that $e\D(\chi)_{W \cdot (v,v^*)} \neq 0$. 
Let $\IC$ be an irreducible component of $\ZC_k^\attractif$ contained 
in the support of $e\D(\chi)$. Then $p \in \IC$ and so 
$p=\lim W \cdot (v,v^*)$. Since $W \cdot (v,v^*) \in (\ZC_k^\t)_\maxi$ 
by the definition of $(\ZC_k^\t)_\maxi$, this implies that $p \in (\ZC_k^\t)_\maxi$, 
as desired.
\end{proof}

\section{Gaudin algebra}\label{sec:gaudin}

\def\gaudin{{\mathrm{Gau}}}
\medskip

\subsection{Definition}
We recall here the definition of {\it Gaudin algebra}~\cite[\S{8.3.B}]{calogero}.
First, let $\CM[V_\reg][W]$ denote the group algebra of $W$ 
over the algebra $\CM[V_\reg]$ (and not the semi-direct 
product $\CM[V_\reg] \rtimes W$). For $y \in V$, let
$$\DC_y^k=\sum_{s \in \Ref(W)} 
\e(s)c_k(s) \frac{\langle y,\a_s\rangle}{\a_s} s\,\,\in \CM[V_\reg][W].$$
Now, let $\gaudin_k(W)$ be the sub-$\CM[V_\reg]$-algebra of 
$\CM[V_\reg][W]$ generated by 
the $\DC_y^k$'s (where $y$ runs over $V$): it will be called the {\it Gaudin algebra} 
({\it with parameter $k$}) associated with $W$. 

Let $\CM(V)$ denote the function field of $V$ (which is the fraction 
field of $\CM[V]$ or of $\CM[V_\reg]$) and let 
$\CM(V) \gaudin_k(W)$ denote the subalgebra $\CM(V) \otimes_{\CM[V_\reg]} \gaudin_k(W)$ 
of the group algebra $\CM(V)[W]$. Recall~\cite[\S{8.3.B}]{calogero} that 
\equat\label{eq:gaudin-com}
\text{\it $\gaudin_k(W)$ is a commutative algebra,}
\endequat
but that $\gaudin_k(W)$ is generally non-split, as shown by the examples 
treated in~\cite[\S{4}]{bonnafe diedral} and~\cite{lacabanne}.

\bigskip

\subsection{Generalized eigenspaces} 
If $v \in V_\reg$, we denote by $\DC_y^{k,v}$ the specialization 
of $\DC_y^k$ at $v$, namely $\DC_y^{k,v}$ is the element of the group 
algebra $\CM W$ equal to
$$\DC_y^k=\sum_{s \in \Ref(W)} 
\e(s)c_k(s) \frac{\langle y,\a_s\rangle}{\langle v, \a_s \rangle} s.$$
Now, if $v^* \in V^*$ and if $M$ is a $\CM W$-module, we define 
$M^{k,v,v^*}$ to be the common generalized eigenspace of the operators 
$\DC_y^{k,v}$ for the eigenvalue $\langle y, v^* \rangle$, for $y$ 
running over $V$. Namely,
$$M^{k,v,v^*}=\{m \in M~|~\forall~y\in V,~
(\DC_y^{k,v} - \langle y,v^* \rangle \Id_M)^{\dim(M)}(m)=0\}.$$
Then 
\equat\label{eq:dec-gaudin}
M = \bigoplus_{v^* \in V^*} M^{k,v,v^*},
\endequat
since $\gaudin_k(W)$ is commutative. 

\bigskip

\begin{lem}\label{lem:fibre-reg}
Let $\chi \in \Irr(W)$ and let $(v,v^*) \in V_\reg \times V^*$. Then 
the following are equivalent:
\begin{itemize}
\itemth{1} $e\D(\chi)_{W \cdot (v,v^*)} \neq 0$.

\itemth{2} $\D^\reg(\chi)_{v,v^*} \neq 0$.

\itemth{3} $E_\chi^{k,v,v^*} \neq 0$. 
\end{itemize}
\end{lem}

\medskip

\begin{proof}
The equivalence between (1) and (2) follows from the Morita equivalence between 
$\CM[V_\reg \times V^*]^W$ and $\CM[V_\reg \times V^*] \rtimes W$ proved 
in~\cite[Lem.~3.1.8(b)]{calogero}. 
Now, as a $\CM[V_\reg]$-module, $\D^\reg(\chi) \simeq \CM[V_\reg] \otimes E_\chi$, 
and the equivalence between~(2) and~(3) follows from the computations 
in~\cite[\S{8.3.B}]{calogero}.
\end{proof}

\medskip

\subsection{Proof of Theorem~A (i.e. Theorem~\ref{theo:chi-tau})}
Let $\chi \in \Irr(W)$ be $\t$-stable and such that $\chit(\t) \neq 0$ 
and let $v \in V_\reg^\t$. By Lemmas~\ref{lem:fibre-fixe} and~\ref{lem:fibre-reg}, 
it is sufficient to show that there exists $v^* \in V^{*\t}$ such that 
$E_\chi^{k,v,v^*} \neq 0$. 

For this, let $\EC$ denote the set of $v^* \in V^*$ such that 
$E_\chi^{k,v,v^*} \neq 0$. Then it follows from~(\ref{eq:dec-gaudin}) that 
$$E_\chi = \bigoplus_{v^* \in \EC} E_\chi^{k,v,v^*}.\leqno{(*)}$$
Since $\t(v)=v$, we have 
$$\lexp{\t}{\DC_y^{k,v}}=
\sum_{s \in \Ref(W)} \e(s) c_k(s) 
\frac{\langle y,\a_s \rangle}{\langle v , \a_s \rangle} \t s \t^{-1} 
= \sum_{s \in \Ref(W)} \e(s) c_k(s) 
\frac{\langle y,\t^{-1}(\a_s) \rangle}{\langle v , \t^{-1}(\a_s) \rangle} s 
= \DC_{\t(y)}^{k,v}.$$
Consequently,
$$\lexp{\t}{E_\chi^{k,v,v^*}}=E_\chi^{k,v,\t(v^*)}.$$
But $\chit(\t)=\Tr(\t,E_\chi) \neq 0$, so 
$\t$ must fix at least one of the generalized eigenspaces in the 
decomposition $(*)$. In other words, 
this implies that there exists $v^* \in \EC$ such 
that $\t(v^*)=v^*$, as desired. The proof is complete.

\bigskip

\section{Complements}

\medskip

\subsection{Conjectures} 
The variety $\ZC_k$ is endowed with a Poisson structure~\cite[\S{1}]{EG} and so the 
variety of fixed points $\ZC_k^\t$ inherits a Poisson structure too, 
as well as all its irreducible components. Recall from Springer Theorem~\ref{theo:springer} 
that $W^\t$ is a reflection group for its action on $V^\t$, so we can define a 
set of pairs $\aleph_\t$ for the pair $(V^\t,W^\t)$ as well as $\aleph$ has been 
defined for the pair $(V,W)$ and, for each parameter $l \in \CM^{\aleph_\t}$, 
we can define a Calogero-Moser space $\ZC_l(V^\t,W^\t)$. 
The following conjecture is a particular 
case of~\cite[Conj.~B]{autosymp} (see~\cite{autosymp} for a discussion 
about the cases where this conjecture is known to hold):

\bigskip

\begin{quotation}
\begin{conj}\label{conj:cm}
Recall that $\t$ is regular. Then there exists a linear map $\l : \CM^\aleph \to \CM^{\aleph_\t}$ 
and, for each $k \in \CM^\aleph$, a $\CM^\times$-equivariant isomorphism of Poisson varieties 
$$\iota_k : (\ZC_k^\t)_\maxi \longiso \ZC_{\l(k)}(V^\t,W^\t).$$
\end{conj}
\end{quotation}

\bigskip

Assume that Conjecture~\ref{conj:cm} holds and keep its notation. Then $\iota_k$ 
restricts to a map 
$\iota_k : (\ZC_k)_\maxi^{\CM^\times} \longiso \ZC_{\l(k)}(V^\t,W^\t)^{\CM^\times}$. 
If $p \in \ZC_{\l(k)}(V^\t,W^\t)^{\CM^\times}$, we denote by $\FG_{\l(k)}^{(\t)}$ 
the corresponding Calogero-Moser $\l(k)$-family of $W^\t$. The next conjecture, still 
inspired by the representation theory of finite reductive groups 
(see again~\cite[Ex.~12.9]{cm-unip} for some explanations), makes 
Conjecture~B more precise:

\bigskip

\begin{quotation}
\begin{conj}\label{conj:jean}
Recall that $\t$ is regular and assume that Conjecture~\ref{conj:cm} holds. 
If $p \in (\ZC_k^\t)_\maxi^{\CM^\times}$, then
$$\sum_{\chi \in \FG_p^\t} |\chit(\t)|^2=\sum_{\psi \in \FG_{\l(k)}^{(\t)}} \psi(1)^2.$$
\end{conj}
\end{quotation}

\bigskip

Note that this last conjecture is compatible with the fact that 
$$\sum_{\chi \in \Irr(W)^\t} |\chit(\t)|^2=|W^\t|=\sum_{\psi \in \Irr(W^\t)} \psi(1)^2,$$
where the first equality follows from the second orthogonality relation for characters.

\bigskip

\subsection{Roots of unity} 
We consider in this subsection a particular (but very important) case 
of the general situation studied in this paper. 
We fix a natural number $d \ge 1$ and a primitive $d$-th root of unity $\z_d$. 
The group of $d$-th roots of unity is denoted by $\mub_d$. An element $w \in W$ 
is called {\it $\z_d$-regular} if the element $\z_d^{-1}w$ of $\Nrm_{\Gb\Lb_\CM(V)}(W)$ 
is regular. In other words, $w$ is $\z_d$-regular if and only if its $\z_d$-eigenspace 
meets $V_\reg$. The existence of a $\z_d$-regular element is not guaranteed: we say 
that $d$ is a {\it regular number} of $W$ if such an element exists. 

\bigskip

\boitegrise{{\bf Hypothesis.} {\it We assume in this subsection, and only in this
subsection, that $d$ is a regular number of $W$. We denote by $w_d$ a 
$\z_d$-regular element and 
we also set $\t_d=\z_d^{-1}w_d$, so that $\t_d$ is a regular element of 
$\Nrm_{\Gb\Lb_\CM(V)}(W)$.}}{0.75\textwidth}

\bigskip

Recall from~\cite{springer} that $w_d$ is uniquely defined up to conjugacy. Note 
that 
\equat\label{eq:vtau}
V^{\t_d}=\Ker(w_d-\z_d\Id_V),\qquad W^{\t_d}=C_W(w_d)\qquad\text{and}\qquad 
\ZC_k^{\t_d}=\ZC_k^{\mub_d}.
\endequat
Since $\t_d$ induces an inner automorphism of $W$, all the irreducible characters 
are $\t_d$-stable. Moreover, if $\chi \in \Irr(W)$, then $\chit(\t_d)=\xi \chi(w_d)$ 
for some root of unity $\xi$, so
$|\chit(\t_d)|^2=|\chi(w_d)|^2$.
This allows to reformulate both Theorem~A and Conjecture~B in this case:

\bigskip

\begin{quotation}
\begin{conj}\label{conj:chi-w}
Recall that $d$ is a regular number. 
Let $p \in \ZC_k^{\CM^\times}$. Then $p$ belongs to $(\ZC_k^{\mub_d})_\maxi$ 
if and only if $\sum_{\chi \in \FG_p} |\chi(w_d)|^2  \neq 0$.
\end{conj}
\end{quotation}

\bigskip

\begin{theo}\label{theo:chi-w}
Recall that $d$ is regular. 
Let $p \in \ZC_k^{\CM^\times}$ be such that 
$\sum_{\chi \in \FG_p^\t} |\chi(w_d)|^2 \neq 0$. 
Then $p$ belongs to $(\ZC_k^{\mub_d})_\maxi$. 
\end{theo}

\bigskip
\def\Part{{\mathrm{Part}}}
\def\Core{{\mathrm{Cor}}}
\def\cor{{\mathrm{cor}}}
\def\quo{{\mathrm{quo}}}
\def\hook{{\mathrm{hk}}}
\def\pard{{\mathrm{par}}_d}

\begin{exemple}[Symmetric group]
We assume here, and only here, that $W=\SG_n$ acting on $V=\CM^n$ by permutation 
of the coordinates, for some $n \ge 2$. The canonical basis of $\CM^n$ is denoted 
by $(y_1,\dots,y_n)$. 
Then there is a unique orbit of hyperplanes, that we denote by $\O$, 
and $e_\O=2$. To avoid too easy cases, we also assume that $k_{\O,0} \neq k_{\O,1}$ 
(so that $\ZC_k$ is smooth~\cite[Cor.~1.14]{EG}) and that $d \ge 2$. 
Saying that $d$ is a regular number is equivalent to say that $d$ divides $n$ or $n-1$. 
Therefore, we will denote by $j$ the unique element of $\{0,1\}$ such that 
$d$ divides $n-j$ and we set $r=(n-j)/d$. 
Then $w_d$ is the product of $r$ disjoint cycles of length $d$, so one can 
choose for instance
$$w_d=(1,2,\dots,d)(d+1,d+2,\dots,2d)\cdots((r-1)d+1,(r-1)d+2,\dots,rd).$$
Then $V^{\t_d}$ is $r$-dimensional, with basis $(v_1,\dots,v_r)$ 
where $v_a=\sum_{b=1}^{d} \z_d^{-b} e_{(a-1)d+b}$ and the group $C_W(w_d) \simeq G(d,1,r)$ 
acting ``naturally'' as a reflection group on $V^{\t_d}=\bigoplus_{a=1}^r \CM v_a$. 

We also need some combinatorics. 
We denote by $\Part(n)$ (resp. $\Part^d(r)$) the set of partitions of $n$ 
(resp. of $d$-partitions of $r$). If $\l \in \Part(n)$, we denote by 
$\chi_\l$ the irreducible character of $\SG_n$ (with the convention of~\cite{geck pfeiffer}: 
for instance $\chi_n=1$ and $\chi_{1^n}=\e$), 
by $\cor_d(\l)$ the $d$-core of $\l$, by $\quo_d(\l)$ 
its $d$-quotient. We let $\Part(n,d)$ denote the set of partitions of $n$ 
whose $d$-core is the unique partition of $j \in \{0,1\}$. Then the map
$$\quo_d : \Part(n,d) \longto \Part^d(r)$$
is bijective. Finally, if $\mu \in \Part^d(r)$, we denote by $\chi_\mu$ the 
associated irreducible character of $C_W(w_d)=G(d,1,r)$ (with the convention of~\cite{geck jacon}). 
It follows from Murnaghan-Nakayama rule that
\equat\label{eq:murnaghan}
\text{$\chi_\l(w_d) \neq 0$ if and only if $\l \in \Part(n,d)$,}
\endequat
and that 
\equat\label{eq:murnaghan-bis}
\chi_\l(w_d)=\pm \chi_{\quo_d(\l)}(1)
\endequat
for all $\l \in \Part(n,d)$ (see for instance~\cite[Page~47]{BMM}).

Now, the smoothness of $\ZC_k$ implies that the map $\ph_k : \Irr(\SG_n) \longto \ZC_k^{\CM^\times}$ 
is bijective (so that Calogero-Moser $k$-families of $\SG_n$ 
are singleton) and it follows from~\cite{bonnafe maksimau} that Conjecture~\ref{conj:cm} holds 
(except that we do not know if the isomorphism respects the Poisson structure), 
so that we have a $\CM^\times$-equivariant isomorphism of varieties
$$\iota_k : (\ZC_k^{\mub_d})_\maxi \longiso \ZC_{\l(k)}(V^{\t_d},G(d,1,r))$$
for some explicit $\l(k) \in \CM^{\aleph_{\t_d}}$. Moreover, $\ZC_{\l(k)}(V^{\t_d},G(d,1,r))$ 
is smooth so that the map $\ph_{\l(k)}^{\t_d} : \Irr(G(d,1,r)) \longto 
\ZC_{\l(k)}(V^{\t_d},G(d,1,r))^{\CM^\times}$ is bijective (that is, Calogero-Moser 
$\l(k)$-families of $G(d,1,r)$ are singleton). Now, by~\cite{bonnafe maksimau}, we have that 
\equat\label{eq:ruslan}
\text{$\ph_k(\chi_\l) \in (\ZC_k^{\mub_d})_\maxi$ if and only if $\l \in \Part(n,d)$,}
\endequat
and that 
\equat\label{eq:ruslan-bis}
\iota_k(\ph_k(\chi_\l))=\ph_{\l(k)}^{\t_d}(\chi_{\quo_d(\l)})
\endequat
for all $\l \in \Part(n,d)$.
Then~\eqref{eq:murnaghan},~\eqref{eq:murnaghan-bis},~\eqref{eq:ruslan} and~\eqref{eq:ruslan-bis} 
show that Conjectures~\ref{conj:chi-w} and~\ref{conj:jean} hold for the symmetric group.\finl
\end{exemple}

\end{document}